\documentclass[12pt,twoside]{article}
\usepackage{amsmath,amssymb,mathrsfs,graphicx}
\usepackage[small,sc]{caption} 
\usepackage[nodvipsnames]{color}

\newtheorem{theorem}{Theorem}[section]
\newtheorem{lemma}[theorem]{Lemma}
\newtheorem{proposition}[theorem]{Proposition}

\newtheorem{definition}[theorem]{Definition}

\newtheorem{rmrk}{Remark} 			

\DeclareMathAlphabet{\mathbfit}{OML}{cmm}{b}{it}

\newenvironment{remark}
{\begin{rmrk} \em}
{\end{rmrk}}

\newcommand{\fn} {function}
\newcommand{\me} {measure}

\newcommand{\erg} {ergodic}
\newcommand{\sy} {system}

\newcommand{\dsy} {dynamical system}

\newcommand{\R} {\mathbb{R}}

\newcommand{\Z} {\mathbb{Z}}
\newcommand{\N} {\mathbb{N}}
\newcommand{\qed} {\hfill {\small Q.E.D.} \par\medskip}
\newcommand{\skippar} {\par\medskip}
\newcommand{\ds} {\displaystyle}
\newcommand{\proof} {\noindent \textsc{Proof.} }
\newcommand{\proofof}[1] {\noindent \textsc{Proof of {#1}.} }
\newcommand{\article}[3] {\textsc{{#1}}, {\itshape {#2}}, {{#3}}.}
\newcommand{\book}[3] {\textsc{{#1}}, {\itshape {#2}}, {{#3}}.}
\newcommand{\vol} {\textbf}
\newcommand{\eps} {\varepsilon}

\newcommand{\rset}[2] {\left\{ #1 \: \left| \: #2 \right. \! \right\} }

\newcommand{\into} {\longrightarrow}

\newcommand{\ind} {\mathcal{J}}   
\newcommand{\leb} {m}   
\newcommand{\sca} {\mathscr{A}}   
\newcommand{\scv} {\mathscr{V}}

\newcommand{\go} {\mathcal{G}}   
\newcommand{\lo} {\mathcal{L}}   

\newcommand{\m} {mixing}
\newcommand{\ob} {observable}

\newcommand{\ui} {(0,1]}   
\newcommand{\rp} {\R_0^+}   
\newcommand{\br} {\tau}   
\newcommand{\ibr} {\phi}   
\newcommand{\avg} {\overline{\mu}}   
\newcommand{\avgleb} {\overline{\leb}}   
\newcommand{\ps} {\mathcal{M}}   

\begin{document}

\title{\textbf{Pomeau-Manneville maps \\ are global-local mixing}}

\author{
\scshape
Claudio Bonanno\thanks{
Dipartimento di Matematica, Universit\`a di Pisa, Largo Bruno 
Pontecorvo 5, 56127 Pisa, Italy. E-mail: 
\texttt{claudio.bonanno@unipi.it}.}
\,and
Marco Lenci\thanks{
Dipartimento di Matematica, Universit\`a di Bologna,
Piazza di Porta San Donato 5, 40126 Bologna, Italy. 
E-mail: \texttt{marco.lenci@unibo.it}.}
\thanks{
Istituto Nazionale di Fisica Nucleare,
Sezione di Bologna, Via Irnerio 46,
40126 Bologna, Italy.}
}

\date{Final version for \\
\emph{Discrete and Continuous Dynamical Systems}\\[5pt]
July 2020}

\maketitle

\begin{abstract}
  We prove that a large class of expanding maps of the unit interval with
  a $C^2$-regular indifferent fixed point in 0 and full increasing branches are
  global-local mixing. This class includes the standard Pomeau-Manneville
  maps $T(x) = x + x^{p+1}$ mod 1 ($p \ge 1$), the Liverani-Saussol-Vaienti 
  maps (with index $p \ge 1$) and many generalizations thereof.

  
\end{abstract}

\section{Introduction}
\label{sec-intro}

By `Pomeau-Manneville map' one generally means a piecewise expanding
map of $[0,1]$ with two increasing surjective branches and an indifferent
fixed point in 0. These maps are so named after Pomeau and Manneville 
who, in the late 1970s, studied numerically approximated versions of them 
to investigate the phenomenon of intermittency in physics \cite{pm, m}.

The first examples of this kind were the maps 
$T_{PM}(x) = x + x^{p+1}$ mod 1, with $p \in \R^+$. 
Throughout this paper we refer to them as \emph{classical 
Pomeau-Manneville (PM) maps}. In 1999, Liverani, Saussol and Vaienti
\cite{lsv} introduced a somewhat simpler one-parameter family:
\begin{equation} \label{lsv}
  T_{LSV}(x) = \left\{ 
  \begin{array}{ll}
    x + 2^p x^{p+1}, & 0 \le x \le \frac12 ; \\[8pt]
    2x-1, & \frac12 < x \le 1,
  \end{array} 
  \right.
\end{equation}
where again $p \in \R^+$ (although in \cite{lsv} only the case $0 < p < 1$ 
was treated). These are nowadays called \emph{LSV maps}. Each of these 
\dsy s has long been known to possess an absolutely continuous invariant 
\me\ $\mu$, unique up to factors, whose density behaves like $x^{-p}$, as 
$x \to 0^+$ \cite{t80}. Thus $\mu$ is infinite if and only if $p\ge 1$.
For such choice of $p$, both the classical PM and the LSV maps have 
become paradigms not just of non-uniformly expanding/hyperbolic maps,
but also of \dsy s preserving an infinite \me. They are undoubtedly the most
common examples in the field of infinite \erg\ theory.

In this paper we prove that a large class of piecewise expanding maps of
the unit interval with an indifferent fixed point in 0 --- including the above
ones and many more --- are global-local \m. This means that for any
$F \in L^\infty(\mu)$ which admits a finite \emph{infinite-volume 
average} 
\begin{equation}
  \avg(F) := \lim_{a \to 0^+} \frac 1 {\mu([a,1])} \int_a^1 F \, d\mu ,
\end{equation}
and any $g \in L^1(\mu)$, one has
\begin{equation}
  \lim_{n \to \infty} \int_0^1 (F \circ T^n) g \, d\mu = \avg(F) \int_0^1 g \, 
  d\mu.
\end{equation}
(In truth, we have a slightly weaker result in the case where the 
\emph{index} of the map $p$ equals 1, see below.) We call all \fn s 
like $F$ `global \ob s' and all \fn s like $g$ `local \ob s'. We 
indicate them, respectively, with uppercase and lowercase letters.

The notion of global-local \m\ has received increasing attention lately
in infinite \erg\ theory; cf.\ \cite{limix, lmmaps, bgl1, bgl2, dn2, ghr} and 
references therein. An example of its usefulness is the fact that, for 
maps with indifferent fixed points, it provides  interesting unconventional
limit theorems, cf.\ \cite[Sect.~3]{bgl1} and \cite[Sect.~3]{bgl2}.

The class of \sy s that we consider here is given by piecewise expanding 
maps $T: \ui \into 
\ui$ with a finite or countable number of increasing surjective branches
and such that, for $x \to 0^+$, $T(x) = x + \kappa x^{p+1} + o(x^{p+1})$,
with $\kappa>0$ and $p \ge 1$.
We also assume a certain condition on the growth of the branches,
see (A5) below, together with standard hypotheses (regularity, distortion
bounds, etc.). The class includes:
\begin{itemize}
\item \emph{generalized classical PM maps}, that is, maps of the form 
$T(x) = x + \kappa x^{p+1}$ mod 1, with $\kappa \in \Z^+$;

\item \emph{generalized LSV maps}, that is, maps whose
branch at 0 is given by $T|_{(0,a_1]}(x) = x + \kappa x^{p+1}$, where
now $\kappa \in \R^+$, and the other branches are 
piecewise linear and increasing;

\item suitable perturbations of the above types.
\end{itemize}

In \cite[Rmks.\ 2.15-2.16]{bgl2} the present authors and P.~Giulietti had 
already discussed global-local \m\ for Pomeau-Manneville maps, albeit 
not in particularly decisive terms.  The problem there was that the proof of 
global-local \m\ required a rather precise knowledge of $h_\mu = 
\frac{d\mu}{d\leb}$, the (infinite) invariant density, which is \emph{not} 
known in general. One of the few cases in which it is known
is that of the map $T(x) = x + x^2$ mod 1, which was therefore proved to 
be global-local \m\ in \cite[Rmk.\ 2.15]{bgl2}. No other case of classical 
PM or LSV map was covered so far.

The techniques used in this work are substantial refinements of those 
employed in \cite{bgl2}, in that they only require to know the behavior of the 
singularity of $h_\mu$ around 0. But a classical result of Thaler (see 
Theorem \ref{thm-thaler}\emph{(a)} below) provides exactly that, for an 
ample class of maps, as a \fn\ of the index $p$. An intermediate step
for our main result (Theorem \ref{main-thm}) is the proof of global-local \m\ 
relative to a \me\ $\nu_p$ whose singular exponent at 0 is exactly 
one unit less than that of $\mu$. What is remarkable is that global-local 
\m\ relative to $\nu_p$ is completely equivalent to global-local \m\ relative 
to $\mu$, for all $p>1$. For $p=1$, the two results are \emph{almost} 
equivalent, in that the class of global \ob s is slightly smaller than the 
optimal class, cf.\ Theorem \ref{main-thm}\emph{(c)}. 

The paper is organized as follows. In Section \ref{sec-setup} we lay out 
the necessary mathematical background and state our results; then
we present a list of examples to which our main theorem applies.
In Section \ref{sec-rp} we prove global-local \m\ for a large class of
maps $\rp \into \rp$ with an ``indifferent fixed point at $\infty$''. Aside
from its own worth, this result is instrumental in the proof of 
the main theorem, which is given in Section \ref{sec-pf}. Finally, 
certain technical results are proved in the Appendix.

\paragraph*{Acknowledgments.} This research was partially supported by 
the PRIN Grant 2017S35EHN of the Ministry of Education, University 
and Research (MIUR), Italy. It is also part of the authors' activity within 
the DinAmicI community \linebreak (\texttt{www.dinamici.org}) and the 
Gruppo Nazionale di Fisica Matematica, INdAM, Italy. C.B.\ was also 
partially supported by the research project 
PRA$\_$2017$\_$22, Universit\`a di Pisa.

\section{Setup and results}
\label{sec-setup}

We now give a precise definition of the class of maps we consider in this
article.
A finite or infinite sequence of numbers $0 = a_0 < a_1 < \ldots < a_k < 
\ldots \le 1$ is given. If the sequence is finite, its last element is denoted 
$a_N$ and it equals $1$; in this case we set $\ind := \{ 0, \ldots, N-1 \}$. 
If the sequence is infinite, $\lim_n a_n = 1$; in this case we set $\ind := \N$. 
For $j \in \ind$, denote $I_j := (a_j, a_{j+1}]$. If $\leb$ denotes the 
Lebesgue \me, then $\{ I_j \}_{j \in \ind}$ is a partition of $\ui$ mod $\leb$,
which acts as the Markov partition of a map $T: \ui \into \ui$ satisfying the 
following conditions:

\begin{itemize} 
\item[(A1)] $T |_{I_j}$ possesses a continuous extension $\br_j : 
  [a_j, a_{j+1}] \into [0,1]$ which is strictly increasing, bijective and 
  twice differentiable on $(a_j, a_{j+1})$. 
  
\item[(A2)] There exist $\kappa > 0$, $p \ge 1$ and $b_o \in (0, a_1)$ 
  such that $\br_0(x) = x + \kappa x^{p+1} + o(x^{p+1})$, as $x \to 0^+$, 
  and $\br_0$ is strictly convex in $[0,b_o]$. This implies in particular 
  that $\br'_0(0) = 1$ and $\br'_0(x) > 1$, for $x \in (0,b_o)$.
  
\item[(A3)] There exists $\Lambda > 1$ such that $T'(x)  \ge
  \Lambda$, for all $x \in [b_o,1) \setminus \{ a_j \}_{j \ge 1}$.
  
\item[(A4)] There exists $K > 0$ such that $\ds \frac{ |T''(x) | }
  { (T'(x))^2 } \le K$, for all $x \in (0,1) \setminus \{ a_j \}_{j \ge 1}$.  
  
\item[(A5)] Set $\ibr_j := \br_j^{-1}$. For all $j \in \ind$, $\ds 
  \sum_{k \ge j} \left( \frac \xi {\ibr_k(\xi)} \right)^{p+1} \!\! \ibr_k'(\xi)$ is 
  a (not necessarily strictly) increasing \fn\ of $\xi \in (0,1)$.
\end{itemize}

We will refer to $\br_j$ as a \emph{branch} of $T$ and to $\ibr_j$ as 
the corresponding \emph{inverse branch}.
Under (much) more general conditions than the above, Thaler 
\cite{t80, t83, t01} proved the following results.

\begin{theorem} \label{thm-thaler} 
  Under the assumptions {\rm (A1)-(A4)}:
  \begin{itemize}
  \item[(a)] $T$ preserves an infinite invariant \me\ $\mu$ which is 
    absolutely continuous w.r.t.\ the Lebesgue \me\ $\leb$, with 
    (infinite) density
    \begin{displaymath}
      h_\mu(x) := \frac{d\mu} {d\leb} (x) = \frac{H_\mu (x)} {x^p},
    \end{displaymath}
    where $H_\mu$ is positive and continuous on $[0,1]$.
  
  \item[(b)] Up to multiplicative constants, $\mu$ is the unique 
    $\leb$-absolutely continuous invariant \me. 
    
  \item[(c)] $T$ is conservative and exact (w.r.t.\ $\leb$ or $\mu$, which 
    is the same).
  \end{itemize}
\end{theorem}

We now introduce two types of \ob s on the space $\ui$. Let $\nu$ be an 
infinite \me\ on $\ui$ such that $\nu([a,1]) < \infty$, for all $a \in \ui$. In 
particular, $\nu$ is $\sigma$-finite. In what follows we refer to any such 
$\nu$ as a \me\ which is \emph{infinite at 0}. Define
\begin{equation} \label{def-go-ui}
  \go(\ui, \nu) := \rset{F \in L^\infty(\ui, \nu) } {\exists \, \overline{\nu} (F) :=
  \lim_{a \to 0^+} \frac1 {\nu([a,1])} \int_a^1 F \, d\nu } .
\end{equation}
We call any $F$ as in the above definition a \emph{global \ob} and say 
that $\overline{\nu} (F)$ is the \emph{infinite-volume average} of 
$F$ relative to $\nu$. We also call any $f \in L^1(\ui, \nu)$ a \emph{local 
\ob}. For added clarity, we denote global, respectively local, \ob s with
uppercase, respectively lowercase, letters. In the remainder we use the 
conventional abbreviation $\nu(f) := \int_0^1 f\, d\nu$.

\begin{definition} \label{def-glm2}
  Given a \me\ $\nu$ which is infinite at 0 and two (sub)classes 
  of global and local
  \ob s, respectively $\go \subseteq \go(\ui, \nu)$ and $\lo \subseteq 
  L^1(\ui, \nu)$, we say that the map $T$ is \textbf{global-local \m}
  w.r.t.\ $\nu, \go, \lo$ if, for all $F \in \go$ and $g \in \lo$,
  \begin{displaymath}
    \lim_{n \to \infty} \nu( (F \circ T^n) g) = \overline{\nu} (F) \nu(g).
  \end{displaymath}
  If the above can be proved for $\go = \go(\ui, \nu)$ and $\lo = 
  L^1(\ui, \nu)$, we say that $T$ is \textbf{fully global-local \m} 
  w.r.t.\ $\nu$.
\end{definition}

\begin{remark} \label{rk-exhaustive}
  In the context of \emph{infinite-volume \m} \cite{limix, lpmu}, 
  there are a few variants of the definition of global-local \m. The one we 
  consider here is the most natural among them and is otherwise 
  denoted \textbf{(GLM2)}. More importantly, in that framework, the class of 
  global \ob s $\go(\ui, \nu)$ corresponds to the \emph{exhaustive family} 
  $\scv := \rset{[a,1]} {0 < a < 1}$. An exhaustive family is a collection of 
  finite-\me\ ``large boxes'' that are used to define the infinite-volume limit. 
  In the present case, where the reference space $\ui$ is one-dimensional 
  and possesses only one ``point at infinity'' (w.r.t.\ $\nu$), $\scv$ is 
  essentially the smallest choice for the exhaustive family, making 
  $\go(\ui, \nu)$ the largest reasonable class of global \ob s, relative to
  $\nu$. Since $L^1(\ui, \nu)$ is the largest class of local \ob s for which
  the definition of global-local \m\ makes sense, this explains the 
  expression `full global-local \m'. 
\end{remark}

\begin{definition} \label{def-ident-avg}
  We say that two \me s $\nu_1$ and $\nu_2$ (that are infinite at 0) 
  give rise to two \textbf{identical} infinite-volume averages $\overline{\nu}_1$ and $\overline{\nu}_2$ if $\go(\ui, \nu_1) = 
  \go(\ui, \nu_2)$ and $\overline{\nu}_1 (F) = \overline{\nu}_2 (F)$ for every 
  $F \in \go(\ui, \nu_1)$. In this case, we write $\overline{\nu}_1 = 
  \overline{\nu}_2$.
\end{definition}

\begin{remark} \label{rk-ident-avg}
  It is obvious that $\overline{\nu}$ only depends on the behavior of $\nu$ 
  around zero, in that, for all $F \in \go(\ui, \nu)$ and $0 < b < 1$, 
  $\overline{\nu}(F) = \overline{\nu}(F 1_{(0,b]})$. In particular, if $\nu_1$ is 
  absolutely continuous w.r.t.\ $\nu_2$ in a neighborhood of 0 and 
  $\frac{d\nu_1}{d\nu_2} (x)$ tends to a positive constant as $x \to 0^+$, 
  then $\overline{\nu}_1 = \overline{\nu}_2$. 
\end{remark}

Definition \ref{def-glm2} of global-local \m\ has a stronger significance, 
in the sense of the decorrelation between a global and a local \ob,
when $\overline{\nu}$ is an invariant functional, i.e., 
$\overline{\nu}(F \circ T^n) = \overline{\nu}(F)$, for all $n \in \N$. This 
is the case, for instance, where $\overline{\nu} = \avg$ and $\mu$ is the 
invariant \me\ guaranteed by Theorem \ref{thm-thaler} \cite[Prop.\ 2.3]{bgl2}.
Of course, we are especially interested in this case. Nevertheless, 
in order to state our main result, Theorem \ref{main-thm} below, we
introduce a special \me\ which will play a role both in the statement and
even more in the proof of the theorem. For $p \ge 1$, let $\nu_p$ be the 
Lebesgue-absolutely continuous \me\ defined by
\begin{equation} \label{h-nu-p}
  h_{\nu_p}(x) := \frac{d\nu_p} {d\leb}(x) = \frac1 {x^{p+1}}
\end{equation}
Clearly, $\nu_p$ is infinite at 0.
\begin{theorem} \label{main-thm}
  Let $T$ satisfy {\rm (A1)-(A5)}. Then:
  \begin{itemize}
  \item[(a)] $T$ is fully global-local mixing relative to $\nu_p$.
  \end{itemize}
  As concerns the invariant \me\ $\mu$:
  \begin{itemize}
  \item[(b)] If $p>1$, $T$ is fully global-local mixing relative to $\mu$ as 
    well.
  \item[(c)] If $p=1$, then $\go (\ui, \nu_1) \subsetneq \go (\ui, \mu)$, with 
    $\avg(F) = \overline{\nu}_1 (F)$ for all $F \in \go (\ui, \nu_1)$. 
    Furthermore, $T$ is global-local mixing w.r.t.\ $\mu$, $\go (\ui, \nu_1)$ 
    and $L^1 (\ui, \mu)$. 
  \end{itemize}
\end{theorem}

To be crystal-clear, the second claim of \emph{(c)} states that, for all $F \in 
\go (\ui, \nu_1)$ and $g \in L^1 (\ui, \mu)$,
\begin{equation}
  \lim_{n \to \infty} \mu((F \circ T^n) g) = \avg(F) \mu(g).
\end{equation}

\begin{remark} \label{main-rk}
Of the five assumptions (A1)-(A5), the last one is certainly the hardest 
to check. A stronger but more readable assumption is:
\begin{itemize} 
\item[(A5)'] For all $j \in \ind$, $\ds \left( \frac {\br_j(x) } x \right)^{p+1} 
  \!\! \frac1 { \br_j'(x) }$ is a (not necessarily strictly) increasing \fn\ of 
  $x \in (a_j, a_{j+1})$.
\end{itemize}
In fact, if (A5)' holds, for all $j$ we can compose the \fn\ 
$( \br_j(x) / x )^{p+1} \! / \br_j'(x)$ with $x = \ibr_j(\xi)$, which is an
increasing \fn\ by (A2). Hence $( \xi / \ibr_j(\xi) )^{p+1} \ibr_j'(\xi)$ 
is increasing in $\xi$ and implies (A5).
\end{remark}

\subsection{Pomeau-Manneville maps and other examples}

We prove Theorem \ref{main-thm} in Section \ref{sec-pf}. For the moment
we verify that several classes of well-known intermittent maps
of the unit interval satisfy its hypotheses.

\paragraph{1. Generalized classical Pomeau-Manneville (PM) maps.}
These are maps of the form 
\begin{equation}
  T(x) = x + \kappa x^{p+1} \mod 1, 
\end{equation}
where $p \ge 1$ and $\kappa \in \Z^+$. The term `generalized' refers 
to the fact that the number $N = \kappa+1$ of branches may be 
larger than one. Assumptions (A1)-(A4) are easily verified. 

We proceed to verify (A5)'.
For all $0 \le j \le N-1$, $a_j$ is the unique $x$ such that $x + 
\kappa x^{p+1} = j$ and $\br_j(x) = x + \kappa x^{p+1} - j$, whence
\begin{equation} \label{ex-10}
  \left( \frac {\br_j(x) } x \right)^{p+1} \!\! \frac1 { \br_j'(x) } = \frac
  {(1 + \kappa x^p - j x^{-1})^{p+1}} {1 + \kappa (p+1) x^p} .
\end{equation}
Verifying (A5)' is tantamount to verifying that the logarithmic derivative 
of the above, w.r.t.\ the variable $z := \kappa x^p$, is non-negative. But
\begin{equation} \label{ex-20}
\begin{split}
  & \frac d {dz} \log \left( \frac {\left(1 + z - j \kappa^{1/p} \, z^{-1/p} 
  \right)^{p+1}} {1 + (p+1) z} \right) \\
  &\qquad = (p+1) \frac{ 1 + j \kappa^{1/p} (1/p) z^{-1-1/p} } 
  {1 + z - j \kappa^{1/p} \, z^{-1/p} } - \frac{p+1} {1 + (p+1)z} \\
  &\qquad = (p+1) \frac{ pz +  j \kappa^{1/p}  z^{-1/p} (z^{-1} + 2p + 1)/p } 
  { \left( 1 + z - j \kappa^{1/p} \, z^{-1/p} \right) (1 + (p+1)z) } > 0,
\end{split}
\end{equation}
for all $z$ such that $1 + z - j \kappa^{1/p} \, z^{-1/p} > 0$, corresponding
to $x > a_j$.

\paragraph{2. Generalized Liverani-Saussol-Vaienti (LSV) maps.} 
These are maps which can have a finite or infinite number of branches. 
The first branch is 
\begin{equation}
  \br_0(x) = x + \kappa x^{p+1}, 
\end{equation}
with $p \ge 1$. Here $\kappa \in \R^+$ with the only obvious constraint 
that $a_1 + \kappa a_1^{p+1} = 1$, for some $0 < a_1 < 1$. The other 
endpoints $a_j$ can be fixed freely and, for $j\ge 1$, 
\begin{equation}
  \br_j(x) = \frac{x - a_j} {a_{j+1} - a_j} , 
\end{equation}
meaning that $\br_j$ is linear, increasing and surjective. (The actual
LSV maps, cf.\ (\ref{lsv}), have $N=2$ branches and are defined by 
$\kappa = 2^p$.)

Again, (A1)-(A4) are easily verified, so we check (A5)'. It is convenient 
to split (A5)' in a number of assumptions, for $j \in \ind$:
\begin{itemize}
  \item[(A5)$'_j$] \ $\ds \left( \frac {\br_j(x) } x \right)^{p+1} \!\! 
  \frac1 { \br_j'(x) }$ is an increasing \fn\ of $x \in (a_j, a_{j+1})$.
\end{itemize}
We observe that (A5)$'_0$ was
already proved in (\ref{ex-10})-(\ref{ex-20}), which in no way depended 
on the fact that $\kappa$ was an integer there. Also, (A5)$'_j$ for  
$j \ge 1$ is immediate since
\begin{equation} 
  \frac {\br_j(x) } x = \frac{1 - a_j x^{-1}} {a_{j+1} - a_j}
\end{equation}
is increasing in $x$ and $\br_j'$ is a positive constant.

\paragraph{3. Perturbations of generalized classical PM maps.} 
With reference to the maps of type 1 above, let $T$ be 
defined by the branches 
\begin{equation}
  \br_j(x) = x + \kappa x^{p+1} -j + \eta_j(x),
\end{equation}
for $j = 0, \ldots, N-1$. Here $\eta_j : [a_j, a_{j+1}] \into \rp$ is twice 
differentiable, $p \ge 1$ as usual and $\kappa \in \R^+$.  We assume 
that  (A1)-(A4) hold. We now prove that, for a sufficiently small $\eps>0$, if
\begin{equation} \label{ex-50}
  | \eta_0(x) | \le \eps x^{2p+1}; \qquad | \eta_0'(x) | \le \eps x^{2p}; \qquad
  | \eta_0''(x) | \le \eps x^{2p-1}, 
\end{equation}
and, for $j\ge 1$,
\begin{equation}  \label{ex-60}
  | \eta_j(x) | \le \eps x^{p+1}; \qquad | \eta_j'(x) | \le \eps x^p; \qquad
  | \eta_j''(x) | \le \eps x^{p-1} ,
\end{equation}
$T$ satisfies (A5)' as well.

Using again the convenient variable $z = \kappa x^p \in [0, \kappa]$, we 
rewrite
\begin{equation} \label{ex-70}
  \left( \frac {\br_j(x) } x \right)^{p+1} \!\! \frac1 { \br_j'(x) } = 
  \frac {\left(1 + z - j \kappa^{1/p} \, z^{-1/p} + u_j(z) \right)^{p+1}} 
  {1 + (p+1) (z + v_j(z) )} ,
\end{equation}
where $u_j(z) := \kappa^{1/p} z^{-1/p} \, \eta_j (\kappa^{-1/p} z^{1/p}) $ and 
$v_j(z) := \eta_j' (\kappa^{-1/p} z^{1/p}) / (p+1)$ (here $\eta_j'$ denotes the
derivative of $\eta_j$ w.r.t.\ $x$). In other words, $\eta_j(x)/x = u_j(z)$ and
$\eta_j'(x) = (p+1) v_j(z)$. The logarithmic derivative of (\ref{ex-70}) 
w.r.t.\ $z$ equals
\begin{equation}
\begin{split}
  & (p+1) \frac{ 1 + j \kappa^{1/p} (1/p) z^{-1-1/p} + u_j'(z) } 
  {1 + z - j \kappa^{1/p} \, z^{-1/p} + u_j(z) } - \frac{(p+1) (1 + v_j'(z))} 
  {1 + (p+1) (z + v_j(z) )} \\
  =: \ & (p+1) \left( \frac{A(z)}{B(z)} - \frac{C(z)}{D(z)} \right) = 
  \frac{p+1} {B(z) D(z)} \left( A(z)D(z) - B(z)C(z) \right),
\end {split}
\end{equation}
where $u_j'$ and $v_j'$ now denote derivatives w.r.t.\ $z$. Observe that,
for the values of $z$ that we are considering, that is, $\kappa a_j^p < z 
< \kappa a_{j+1}^p$, both $B(z)$ and $D(z)$ are positive, as
they correspond, respectively, to $\br_j(x)/x$ and $\br_j'(x)$, with 
$a_j < x < a_{j+1}$. Proving that $A(z)D(z) - B(z)C(z) \ge 0$, for the 
mentioned values of $z$, will thus establish (A5)$'_j$. After some computations 
and regrouping of similar terms, one has:
\begin{equation} \label{ex-80}
\begin{split}
  & A(z)D(z) - B(z)C(z)  \\
  &\quad = \left[ u_j'(z) + (p+1) v_j(z) - u_j(z) - v_j'(z) + (p+1) u_j'(z) v_j(z) - 
  u_j(z) v_j'(z) \right] \\
  &\qquad + z  \left[ p + (p+1) u_j'(z) - v_j'(z) \right] \\
  &\qquad + z^{-1/p} j \kappa^{1/p} \left[ 2 + 1/p + v_j'(z) \right] \\
  &\qquad + z^{-1-1/p} j \kappa^{1/p} \left[ 1/p + (1 + 1/p) v_j(z) \right] .
\end {split}
\end{equation}

Let us verify (A5)$'_0$. In terms of $u_0$ and $v_0$, the conditions (\ref{ex-50}) 
read:
\begin{equation} 
  | u_0(z) | \le \eps_o z^2; \quad | u_0'(z) | \le \eps_o z; \quad
  | v_0(z) | \le \eps_o z^2, \quad | v_0'(z) | \le \eps_o z,
\end{equation}
for some $\eps_o > 0$ proportional to $\eps$. Here we have used that
$\frac{dx}{dz} = \frac1p \kappa^{1/p} z^{(1-p) / p}$. So, for $\eps$ (and thus $\eps_o$) 
small enough, the contribution of all the terms in the second and third lines  
of (\ref{ex-80}) can be made bigger than, say, $(p/2)z > 0$, for all $z \in
(0, \kappa \, a_1^{p})$. Since the fourth and fifth lines of 
(\ref{ex-80}) are null for $j=0$, this gives (A5)$'_0$.

In the case $j \ge 1$, the conditions (\ref{ex-60}) imply that
\begin{equation} 
  | u_j(z) | \le \eps_o z; \quad | u_j'(z) | \le \eps_o; \quad
  | v_j(z) | \le \eps_o z, \quad | v_j'(z) | \le \eps_o,
\end{equation}
for a small $\eps_o$. In this case we must consider values of $z \in
(\kappa a_j^p, \kappa a_{j+1}^p)$. Since in this case $pz > p \kappa a_j^p 
> 0$, the above hypotheses are 
enough to ensure that the second and third lines of (\ref{ex-80}) can be 
made bigger than a positive constant. By the same reasoning, and the fact
that $z$ is bounded above, we can say the same about the fourth and fifth 
lines. This proves (A5)$'_j$.

\paragraph{4. Perturbations of generalized LSV maps.}
In view of the maps of type 2 above, consider a $T$ given by 
\begin{align}
  \br_0(x) &= x + \kappa x^{p+1} + \eta_0(x) ; \\
  \br_j(x) &= \frac{x - a_j} {a_{j+1} - a_j} + \eta_j(x) \qquad (j \ge 1).
  \label{brj-4}
\end{align}
Once again, we assume (A1)-(A4). (In particular, we have $\eta_j(a_j) = 
\eta_j(a_{j+1}) = 0$.) In addition, suppose that $\eta_0$ verifies conditions 
(\ref{ex-50}), with the same $\eps$ determined \emph{a fortiori} to work
for the example 3, and, for $j \ge 1$,
\begin{equation} \label{last-conds}
  \eta_j''(x) \le 0; \qquad \eta_j(x) - x \eta_j'(x) \le \frac{a_j} {a_{j+1} - a_j} .
\end{equation}
Let us observe that it must be $0 \le \eta_j(x) \le 1$, the first inequality coming
from the concavity of $\eta_j$ and the second from (\ref{brj-4}) and (A1).

The above hypotheses imply (A5)'. In fact, (A5)$'_0$ holds by construction. 
As for (A5)$'_j$, with $j\ge 1$, observe that $\br_j'' = \eta_j'' \le 0$, so that 
$1/ \br'_j$ is increasing. Furthermore,
\begin{equation} 
  \frac d {dx} \left( \frac{\br_j(x)} x \right) = \frac{a_j} {a_{j+1} - a_j} \, 
  \frac1 {x^2} + \frac{\eta_j'(x)} x - \frac{\eta_j(x)} {x^2}
\end{equation}
is non-negative if and only if the second inequality of (\ref{last-conds})
holds.

\bigskip

We conclude this section by emphasizing that in no way are the conditions
presented in examples 3 and 4 necessary for (A5)'. They have been chosen 
only because of their simplicity and capacity to generate many examples.
When faced with a specific map $T$, the most reasonable course of action 
is to simply check (A1)-(A5) or, if calculating the inverse branches of $T$ 
is too cumbersome, (A1)-(A4) and (A5)'.

\section{Maps on the half-line}
\label{sec-rp}

In this section we introduce a class of Markov maps of the half-line 
$\rp := [0,+\infty)$ that are analogous to the maps of the unit interval 
discussed earlier. For these maps we prove full global-local \m\ w.r.t.\ a 
class of \me s, including in many cases the Lebesgue \me. This result is 
interesting in its own right and will be the basis for the proof of Theorem 
\ref{main-thm}, which we present in the next section.

A map $T: \rp \into \rp$ is defined as follows. There exists 
a finite or infinite sequence of numbers $a_1 > a_2 > \ldots > 
a_k  > \ldots \ge 0$. If the sequence is finite, its last element is 
$a_N := 0$; in this case we set $\ind := \{ 0, \ldots, N-1 \}$. If the 
sequence is infinite, $\lim_n a_n = 0$; in this case we set 
$\ind := \N$. For $j \in \ind$, denote $I_j := [a_{j+1}, a_j)$,
where we have conventionally put $a_0 := +\infty$. The following 
assumptions hold:

\begin{itemize} 
\item[(B1)] For all $j \in \ind$, $ \br_j := T |_{I_j}$ is an 
increasing diffeomorphism $I_j \into \rp$, up to $a_{j+1}$ and 
$0$, on the domain and codomain, respectively. 

\item[(B2)] $T$ is exact w.r.t.\ $\leb$, the Lebesgue \me\ on
$\rp$.

\item[(B3)] Set $\ibr_j := \br_j^{-1}$. For all $j \in \ind$, $\ds \sum_{k \ge j} 
\ibr_k'$ is a (not necessarily strictly) decreasing \fn\ of $\rp$.
\end{itemize}

Let $\nu$ be an infinite, locally finite, \me\ on $\rp$. In analogy to 
(\ref{def-go-ui}), we define the class of global \ob s relative to $\nu$ to be
\begin{equation} \label{def-go-rp}
  \go(\rp, \nu) := \rset{F \in L^\infty(\rp, \nu) } {\exists \, \overline{\nu} (F) :=
  \lim_{a \to +\infty} \frac1 {\nu([0,a])} \int_0^a F \, d\nu } .
\end{equation}
Again we call local \ob\ any $f \in L^1(\rp, \nu)$, and use the abbreviation
$\nu(f) = \int_0^\infty f \, d\nu$. In view of Remark \ref{rk-exhaustive},
we emphasize that the above choice of global \ob s corresponds to the 
exhaustive family $\scv := \rset{[0,a]} {a > 0}$ for $\rp$, and it is effectively
the largest class of global \ob s relative to $\nu$. The definitions of 
global-local \m\ and identical infinite-volume averages (Definitions 
\ref{def-glm2} and \ref{def-ident-avg}) naturally carry over to this setting.

In this context $\leb$ denotes the Lebesgue \me\ on $\rp$. Also, for $0 < q 
\le 1$, we introduce the Lebesgue-absolutely continuous \me\ $\lambda_q$ 
defined by the (infinite) density
\begin{equation} \label{def-lambda-q}
  h_{\lambda_q}(y) := \frac{d\lambda_q} {d\leb} (y) = \frac 1 {(1+y)^q}
\end{equation}
The rest of the section is devoted to proving the following result:

\begin{theorem} \label{thm-rp}
  Under assumptions {\rm (B1)-(B3)}:
  \begin{itemize}
  \item[(a)] $T$ is fully global-local mixing relative to $\leb$.
  \item[(b)] $T$ is fully global-local mixing relative to $\lambda_q$, for all
    $q \in (0,1)$.
  \item[(c)] $\go(\rp, \leb) \subsetneq \go(\rp, \lambda_1)$, with
    $\overline{\lambda}_1(F) = \avgleb(F)$ for all $F \in \go(\rp, \leb)$. 
    Moreover, $T$ is global-local \m\ w.r.t.\ $\lambda_1$, $\go(\rp, \leb)$, 
    $L^1(\rp, \lambda_1)$.
  \end{itemize}
\end{theorem}  

\proof Let us start with assertion \emph{(a)}. We use the same technique 
as in the proof of \cite[Thm.~5.2]{bgl2}. We still write all the arguments 
because the proof in \cite{bgl2} was given for the case where $T$ is a map 
of the unit interval and the \me\ used there was invariant, and not just 
non-singular.

In what follows we write $L^1$ for $L^1(\rp, \leb)$. Let $P = P_{T,\leb}: 
L^1 \into L^1$ be the Perron-Frobenius operator of $T$, i.e., the transfer 
operator of $T$ relative to the Lebesgue \me\ $\leb$. It is well known that 
$P$ acts as
\begin{equation} \label{pf-op}
  P g(y) = \sum_{j \in \ind} \ibr_j'(y) \, g( \ibr_j(y)).
\end{equation}
Also, it is a positive operator and, for all $g \ge 0$, $\| P g \|_1 = \| g \|_1$. 
Thus, for a general $g \in L^1$, $\| P g \|_1 \le \| g \|_1$. In our case
$P$ enjoys this crucial property as well:

\begin{lemma} \label{main-lem}
  If $g \in L^1$ is a decreasing (thus necessarily non-negative) \fn\ of $\rp$, 
  then so is $P g$. In other words, $P$ preserves the cone of the 
  decreasing \fn s in $L^1$.
\end{lemma}

\proofof{Lemma \ref{main-lem}} The above-mentioned cone is spanned 
by the finite linear combinations, with positive coefficients, of the \fn s
$1_{[0,a]}$, for $a > 0$. Since $P$ is a contraction operator, namely 
$\| P\| \le 1$, it suffices to prove the assertion for $g = 1_{[0,a]}$. 

Let $j \in \ind$ be the unique index such that $a \in I_j$ and set
$b := \br_j(a)$. By (\ref{pf-op})
\begin{equation} 
  P 1_{[0,a]} (y) = \left\{ \begin{array}{ll}
    \ds \sum_{k \ge j} \ibr_k' (y), & y < b ; \\[18pt]
    \ds \sum_{k \ge j+1} \ibr_k' (y), & y \ge b. 
  \end{array} \right.
\end{equation}
Now, the two branches on the above r.h.s.\ are decreasing \fn s of $y$
by (B3) and the gap between them, at $y=b$, is $-\ibr_j' (b)$, which is
negative by (B1).
\qed

Lemma \ref{main-lem} guarantees that every decreasing $g \in L^1$ is 
a \emph{persistently decreasing} local \ob, relative to $P$. This means that, 
for all $n\in\N$, $P^n g$ is decreasing. (The notion of \emph{persistently
monotonic} \ob\ was introduced in \cite{bgl2}.) 

We now need the following general lemma:

\begin{lemma} \label{lem-pmu-llm}
  Let $T$ be a non-singular, exact  endomorphism of the $\sigma$-finite, 
  infinite \me\ space $(\ps, \sca, \nu)$. Then:
  \begin{itemize}
  
  \item[(a)] Given $F \in L^\infty$, if there exists $\overline{F} \in \R$ such that 
  the limit
  \begin{displaymath}
    \lim_{n \to \infty} \nu( (F \circ T^n) g) = \overline{F} \nu(g)  
  \end{displaymath}
  holds for some $g \in L^1$, with $\nu(g) \ne 0$, then it holds for
  all $g \in L^1$.
  
  \item[(b)] $T$ is \textbf{local-local \m}, that is, for all $f \in L^\infty \cap L^1$ 
    and $g \in L^1$, 
    \begin{displaymath}
      \lim_{n \to \infty} \nu((f \circ T^n) g) = 0.
    \end{displaymath}
  \end{itemize}
\end{lemma}

The above statements were essentially proved in \cite[Lem.~3.6 and Thm.~3.5(b)]{lpmu}. However, since they were stated with stronger 
hypotheses there, and also in the interest of completeness, we present the
proof of Lemma \ref{lem-pmu-llm} in the Appendix \ref{app-lem-pmu-llm}.

\skippar

Assumption (B2) and Lemma \ref{lem-pmu-llm}\emph{(a)} show that it 
suffices to verify global-local \m\ for a single persistently decreasing local 
\ob\ $g$ with $\| g \|_1 = 1$. Also, by possibly centering the global \ob\ (i.e., 
using $F - \avgleb(F)$ instead of $F$), we can always assume that 
$\avgleb(F) = 0$. Since, by definition of $P$, $\leb((F\circ T^n)g) = 
\leb(F \, P^n g)$, it remains to verify that
\begin{equation} \label{to-prove}
  \lim_{n \to \infty} \leb(F \, P^n g) = 0.
\end{equation}

Fix $\eps > 0$. By definition of $\avgleb$ we can find $M \in \R^+$ such that,
for all $a \ge M$,
\begin{equation} \label{main20}
  \frac 1 a \left| \int_0^a F(y) \, dy \right| < \frac{\eps}2.
\end{equation}
For $y \in \rp$ and $n \in \N$, set
\begin{equation} \label{main30}
  \gamma_n(y) = \gamma_{n,M}(y) := \min\{ P^n g(M), P^n g(y) \}.
\end{equation}
Since $g$ is persistently decreasing, $\gamma_n$ is a positive, 
(non strictly) decreasing \fn, which is constant on $[0,M]$. It is a local
\ob\ because $\| \gamma_n \|_1 \le \| P^n g \|_1 =  \| g \|_1 = 1$. We 
have
\begin{equation} \label{main40}
  \leb( F \, P^n g) = \int_0^\infty F \gamma_n \, d\leb +
  \int_0^M F (P^n g - \gamma_n) \, d\leb =: \mathcal{I}_1 + 
  \mathcal{I}_2.
\end{equation}

To estimate $\mathcal{I}_2$, let us notice that
\begin{equation} 
  0 \le \int_0^M (P^n g - \gamma_n) \, d\leb \le \int_0^M P^n g \, d\leb 
  = \leb \! \left( (1_{[0,M]} \circ T^n) g \right)
\end{equation}
As $n \to \infty$, this term vanishes by Lemma \ref{lem-pmu-llm}\emph{(b)}.
So, for all $n$ large enough,
\begin{equation} \label{main50}
  | \mathcal{I}_2 |  \le \| F \|_\infty  \int_0^M (P^n g - \gamma_n) \, d\leb \le 
  \frac{\eps}2.
\end{equation}
Let us now study $\mathcal{I}_1$. We introduce the generalized inverse 
of the \fn\ $\gamma_n$, which we define to be $\gamma_n^{-1}(r) := \inf 
\rset{y \in \rp} {\gamma_n(y) \le r}$. Evidently, $\gamma_n^{-1}$ is 
decreasing and $\gamma_n^{-1}(r)=0$ for $r \ge \gamma_n(M)$. We 
rewrite $\mathcal{I}_1$ as a double integral:
\begin{equation}
  \mathcal{I}_1 = \int_0^\infty F(y) \left( \int_0^{\gamma_n(y)} \!\!\ dr 
  \right) dy = \int_0^{\gamma_n(M)} \left( \int_0^{\gamma_n^{-1}(r)} F(y) \, 
  dy \right) dr ,
\end{equation} 
where in the second equality we have used Fubini's Theorem to 
interchange the order of integration. Therefore, using (\ref{main20}) with 
$a := \gamma_n^{-1}(r) \ge M$, we obtain
\begin{equation} \label{main80}
\begin{split}
  | \mathcal{I}_1 | &\le \int_0^{\gamma_n(M)} \left| 
  \int_0^{\gamma_n^{-1}(r)} F(y) \, dy \right| dr \\
  & \le \frac{\eps}2 \int_0^{\gamma_n(M)} 
  \int_0^{\gamma_n^{-1}(r)} dy \ dr \\
  &= \frac{\eps}2 \int_0^\infty \int_0^{\gamma_n(y)} \!\!\ dr \ dy \\
  &= \frac{\eps}2 \leb(\gamma_n) \le \frac{\eps}2.
\end{split}
\end{equation} 
Observe that this estimate holds for all $n$. Together with 
(\ref{main50}) and (\ref{main40}), it proves (\ref{to-prove}) and
thus assertion \emph{(a)} of the theorem.

\skippar

The other assertions will follow easily from the following
lemma, which is proved in Appendix \ref{app-claudio}.

\begin{lemma} \label{lem-claudio}
  In view of definitions (\ref{def-go-rp})-(\ref{def-lambda-q}):
  \begin{itemize}
  \item[(a)] For $q \in (0,1)$, $\go(\rp, \lambda_q) = \go(\rp, \leb)$ and 
    $\overline{\lambda}_q(F) = \avgleb(F)$, for all $F \in \go(\rp, \lambda_q)$.
    
  \item[(b)] $\go(\rp, \leb)\subsetneq \go(\rp, \lambda_1)$ and 
    $\overline{\lambda}_1(F) = \avgleb(F)$, for all $F \in \go(\rp, \leb)$.
  \end{itemize}
\end{lemma}

In fact, in all cases $0 < q \le 1$, take any $F \in \go(\rp, \leb) \subseteq 
\go(\rp, \lambda_q)$ and $g \in L^1(\rp, \lambda_q)$. By Theorem 
\ref{thm-rp}\emph{(a)} and Lemma \ref{lem-claudio},
\begin{equation}
  \lim_{n \to \infty} \lambda_q ((F \circ T^n) g) = \lim_{n \to \infty} \leb 
  ((F \circ T^n) g h_{\lambda_q}) = \avgleb(F) \leb ( g h_{\lambda_q} ) = 
  \overline{\lambda}_q(F) \lambda_q (g). 
\end{equation}
This simultaneously gives assertions \emph{(b)} and \emph{(c)}
of Theorem \ref{thm-rp}.
\qed

\section{Proof of the main theorem}
\label{sec-pf}

We have already mentioned that the maps studied in Section \ref{sec-rp}
are analogous to the maps on the unit interval that are the object 
of this proof. Our strategy, in fact, will be to reduce the problem to an
application of Theorem \ref{thm-rp}, by means of a suitable conjugation. 
Of the infinitely many isomorphisms between $\ui$ and $\rp$, it turns
out that one of the most convenient for our purposes is $\Psi : \ui \into \rp$, 
defined by
\begin{equation} \label{psi}
  \Psi(x) := \int_x^1 h_{\nu_p} (\xi) \, d\xi =
  \int_x^1 \xi^{-p-1} \, d\xi = \frac{x^{-p}-1} p .
\end{equation}

Throughout this section we indicate with a subscript $_o$ all objects 
pertaining to the space $\rp$. Specifically:
\begin{enumerate}
\item Given a  map $T: \ui \into \ui$, we denote by $T_o := 
  \Psi \circ T \circ \Psi^{-1}$ its conjugate on $\rp$.
  
\item Given a \me\ $\nu$ on $\ui$ that is infinite at 0, its push-forward 
  $\nu_o := \Psi_* \nu = \nu \circ \Psi^{-1}$ is an infinite, locally finite \me\ 
  on $\rp$. Observe that if $\nu$ is Lebesgue-absolutely continuous, then 
  so is $\nu_o$.
  
\item For any local \ob\ $f \in L^1(\ui, \nu)$ it follows from the definition of
  $\nu_o$ that $f_o := f \circ \Psi^{-1} \in L^1(\rp, \nu_o)$. 

\item For any global \ob\ $F \in \go(\ui, \nu)$ it is readily verified that the 
  corresponding \ob\ $F_o := F \circ \Psi^{-1}$ belongs to $\go(\rp, \nu_o)$, 
  and $\overline{\nu}(F) = \overline{\nu}_o(F_o)$. 
\end{enumerate}
The map $f \mapsto f_o$ defined in points 3 and 4 for local and
global \ob s, respectively, is clearly a bijection. Let us call $\lo_o$ and
$\go_o$, respectively, the images of $\lo$ and $\go$. This shows that
the isomorphism of non-singular \dsy s $(\ui, \nu, T) \cong (\rp, \nu_o, T_o)$
also preserves the infinite-volume structure and the property of  
global-local \m. More precisely, the former \sy\ is global-local \m\ relative
to $\nu, \go, \lo$ if and only if the latter is global-local \m\ relative
to $\nu_o, \go_o, \lo_o$. Hence, the former is fully global-local \m\ 
relative to $\nu$ if and only if the latter is fully global-local \m\ relative
to $\nu_o$. 

Let us therefore prove that, for any $T$ satisfying the hypotheses of 
Theorem \ref{main-thm}, the conjugate $T_o$ verifies (B1)-(B3). (B1)
is immediate since $\Psi$ is a diffeomorphism. (B2) comes from
Theorem \ref{thm-thaler}. As for (B3), we observe that the inverse 
branches of $T_o$ are related to those of $T$ by the same conjugation,
i.e., $\ibr_{o,k} = \Psi \circ \ibr_k \circ \Psi^{-1}$, whence
\begin{equation}
\begin{split}
  \ibr_{o,k}'(y) &= \Psi' \!\left( \ibr_k( \Psi^{-1}(y)) \right) \, \ibr_k' \!\left( 
  \Psi^{-1}(y) \right) \, \frac1 {\Psi' \!\left( \Psi^{-1} (y) \right) } = \\
  &= \left( \frac{\ibr_k( \Psi^{-1}(y))} {\Psi^{-1} (y)} \right)^{-p-1} \!
  \ibr_k' \!\left( \Psi^{-1}(y) \right) .
\end{split}
\end{equation}
Here we have used that $\Psi'(x) = -x^{-p-1}$. Summing over $k \ge j$ and 
composing with $y = \Psi(\xi)$, which is a decreasing \fn\ of $\xi$, proves 
that (B3) and (A5) are equivalent statements. Thus Theorem \ref{thm-rp} 
applies. This implies that $T$ is fully global-local \m\ relative to both
$\Psi^{-1}_* \leb_o$ and $\Psi^{-1}_* \lambda_q$, for all $q \in (0,1)$,
and it is global-local \m\ w.r.t.\ $\Psi^{-1}_* \lambda_1$, 
$\go(\ui, \Psi^{-1}_* \leb_o)$, $L^1(\ui, \Psi^{-1}_* \lambda_1)$.

Now, recalling the definition (\ref{h-nu-p}) of $\nu_p$, a straightforward 
calculation based on (\ref{psi}) shows that, for any Lebesgue-absolutely 
continuous \me\ $\nu_o$ on $\rp$, 
\begin{equation}
  \frac{d(\Psi^{-1}_* \nu_o)} {d\nu_p} = \frac{d\nu_o} {d\leb_o} \circ \Psi,
\end{equation}
where $\leb_o$ is the Lebesgue \me\ on $\rp$. This implies in particular 
that 
\begin{equation} \label{final-20}
  \Psi^{-1}_* \leb_o = \nu_p. 
\end  {equation}
Also, in light of (\ref{def-lambda-q}), for all $q \in (0,1]$ and $x \in \ui$ we 
have
\begin{equation}
  \frac{d(\Psi^{-1}_* \lambda_q)} {d\nu_p}(x) = \frac 1 {(1 + \Psi(x))^q}.
\end{equation}
As $x \to 0^+$, the above expression is asymptotic to $(px^p)^q$,
cf.\ (\ref{psi}). Thus, for $x \to 0^+$, 
\begin{equation}
  \frac{d(\Psi^{-1}_* \lambda_q)} {d\leb}(x) \sim p^q \, x^{pq-p-1},
\end{equation}
whence, by Theorem \ref{thm-thaler}\emph{(a)},
\begin{equation}
  \frac{d(\Psi^{-1}_* \lambda_{1/p})} {d\mu}(x) \sim \frac{p^{1/p}} {H_\mu(0)}.
\end{equation}

Remark \ref{rk-ident-avg} then shows that $\overline{ 
\Psi^{-1}_* \lambda_{1/p} } = \avg$, in the sense of Definition
\ref{def-ident-avg}. This fact and (\ref{final-20}) show that statements 
\emph{(a), (b), (c)} of Theorem \ref{main-thm} come from the corresponding 
statements of Theorem \ref{thm-rp}, with $q = 1/p$.
\qed

\appendix

\section{Appendix: Technical results}

\subsection{Proof of Lemma \ref{lem-pmu-llm}}
\label{app-lem-pmu-llm}

We start by recalling a famous result by Lin \cite{li}: $T$ is exact
if and only if $\| P^n f \|_1 \to 0$, as $n \to \infty$, for all $f \in L^1$ with 
$\nu(f) = 0$. Thus, for all $F \in L^\infty$ and $g \in L^1$ with $\nu(g)=0$, 
\begin{equation} \label{glm1}
  \lim_{n \to \infty} \left| \nu((F \circ T^n) g) \right| = \lim_{n \to \infty} 
  \left| \nu(F \, P^n g) \right| \le \| F \|_\infty \lim_{n \to \infty}  \| P^n g \|_1 = 0.
\end{equation}
The property that $\nu((F \circ T^n) g)$ vanishes, as $n \to \infty$, for all 
$F \in \go$ and $g \in \lo$ is called \textbf{(GLM1)} (w.r.t.\ $\nu, \go, \lo$); 
cf.\ \cite{lpmu, lmmaps, bgl2}.

\begin{lemma}
  \label{lem-pmu}
  Under the hypotheses of Lemma \ref{lem-pmu-llm}, consider $F \in L^\infty$. 
  If, for some $\ell \in \R$ and $\eps \ge 0$, the limit
  \begin{displaymath}
    \limsup_{n \to \infty} \left| \frac{ \nu((F \circ T^n) g) } {\nu(g)} 
    - \ell \right| \le \eps
  \end{displaymath}
  holds for some $g \in L^1$ (with $\nu(g) \ne 0$), then it holds for
  all $g \in L^1$ (with $\nu(g) \ne 0$).
\end{lemma}

\proofof{Lemma \ref{lem-pmu}} Suppose the above limit holds for $g_0
\in L^1$. Take any other $g \in \lo$, with $\nu(g) \ne 0$. We have:
\begin{equation}
\begin{split}
  & \left| \frac{ \nu((F \circ T^n) g) } { \nu(g) } - \ell \right| \\
  &\qquad \le \left| \nu \! \left( (F \circ T^n) \left( \frac{g}{\nu(g)} -
  \frac{g_0}{\nu(g_0)} \right) \right) \right| + \left| 
  \frac{ \nu((F \circ T^n) g_0) } { \nu(g_0) } - \ell \right| .
\end{split}
\end{equation}
The first term of the above r.h.s.\ vanishes, as $n \to \infty$, due to 
(\ref{glm1}).
\qed

At this point assertion \emph{(a)} of Lemma \ref{lem-pmu-llm} is all but 
proved. For any $g \in L^1$ with $\nu(g) \ne 0$ one applies the above 
lemma with $\ell := \overline{F}$ and $\eps := 0$. In the case $\nu(g) = 0$, 
one applies (\ref{glm1}) directly.

As for \emph{(b)}, observe that, since $(\ps, \sca, \nu)$ is $\sigma$-finite 
and infinite, one can always find a set $A \in \sca$ whose \me\ is finite but 
larger than any predetermined number. Thus, for any $\eps > 0$, there
exists $A \in \sca$ with $\| f \|_1 / \eps < \nu(A) < \infty$. Set $g_\eps
:= 1_A / \nu(A)$. We have that
\begin{equation}
  \left| \frac{ \nu( (f \circ T^n) g_\eps) } {\nu(g_\eps)} \right|  = 
  | \nu( (f \circ T^n) g_\eps) | \le \| f \|_1 \, \| g_\eps \|_\infty 
  \le \eps.
\end{equation}
By Lemma \ref{lem-pmu},
\begin{equation}
  \limsup_{n \to \infty} \left| \frac{ \nu( (f \circ T^n) g) } {\nu(g)}
  \right| \le \eps
\end{equation}
holds for all $g \in L^1$ with $\nu(g) \ne 0$.  Since $\eps$ is arbitrary, we 
get that the above r.h.s.\ is zero. When $\nu(g) = 0$, the fact that 
$\lim_{n \to \infty} \nu( (f \circ T^n) g) = 0$ follows directly from (\ref{glm1}). 
This ends the proof of Lemma \ref{lem-pmu-llm}.
\qed

\subsection{Proof of Lemma \ref{lem-claudio}}
\label{app-claudio}

We start with some preliminary results.

\begin{proposition}\label{prop-1-app-c}
  For $q \in (0,1]$, $\go(\rp,\leb) \subseteq \go(\rp,\lambda_q)$ 
  and $\overline{\lambda}_q(F) = \avgleb(F)$, for all $F \in \go(\rp, \leb)$.
\end{proposition}

\proof Let $F\in \go(\rp,\leb)$. Without loss of generality we suppose 
$\avgleb(F)=0$. Therefore, it is enough to show that 
\begin{equation}\label {voglio-q}
  \lim_{a\to \infty}\, \frac{1}{\lambda_q([0,a])} \int_0^a\, F\, d\lambda_q = 0.
\end{equation}
Using integration by parts we have
\begin{equation} \label{serve-dopo}
\begin{split}
   \int_0^a\, F\, d\lambda_q &= \int_0^a\, F(y)\, \frac{1}{(1+y)^q}\, dy \\
  &= (1+a)^{-q} \int_0^a F(y)\, dy + q \int_0^a \frac{1}{(1+y)^{q+1}}\, 
  \Big(\int_0^y\, F(s)\, ds\Big)\, dy .
\end{split}
\end{equation}
Now fix any $\eps>0$. By definition of $\avgleb$, there exists $M=M(\eps) 
\in \R^+$ such that, for all $a > M$, 
\begin{equation} \label{iv-leb-1}
  \left| \int_0^a F(y)\, dy \right| < \eps a.
\end{equation}
For all  $a > M$ we write
\begin{equation} \label{interm-q}
\begin{split}
  \frac{1}{\lambda_q([0,a])} \int_0^a F\, d\lambda_q &= \frac{(1+a)^{-q}}
  {\lambda_q([0,a])} \int_0^a F(y)\, dy \\
  &\quad + \frac{q}{\lambda_q([0,a])} \int_0^{M} \frac{(\int_0^y F(s)\, ds)}
  {(1+y)^{q+1}}\, dy \\
  &\quad + \frac{q}{\lambda_q([0,a])} \int_{M}^a \frac{(\int_0^y F(s)\, ds)}
  {(1+y)^{q+1}}\, dy .
\end{split}
\end{equation}

Let us now specialize to the case $q \in (0,1)$. By (\ref{iv-leb-1}) and 
using that 
\begin{equation} \label{lambdaq-0a}
  \lambda_q([0,a]) = \frac{(1+a)^{1-q}-1} {1-q} ,
\end{equation}
we get
\begin{equation} 
\begin{split}
  \left| \frac{1}{\lambda_q([0,a])} \int_0^a F\, d\lambda_q \right| &\le
   \eps\,  \frac{a (1+a)^{-q}}{\lambda_q([0,a])} \\
   &\quad + \frac{q}{\lambda_q([0,a])} \left| \int_0^{M} \frac{(\int_0^y F(s)\, ds)}
   {(1+y)^{q+1}}\,  dy \right| \\
   &\quad  + \frac{q}{\lambda_q([0,a])} \int_{M}^a \frac{\eps y}{(1+y)^{q+1}}\, dy \\
   &\le \eps\, \frac{(1-q) a}{1+a-(1+a)^q} \\
   &\quad + \frac{q}{\lambda_q([0,a])} \left| \int_0^{M} \frac{(\int_0^y F(s)\, ds)}
   {(1+y)^{q+1}}\, dy \right| \\
   &\quad + \eps\, \frac{q}{\lambda_q([0,a])} \int_{M}^a \frac{1}{(1+y)^q}\, dy \\
   &= \eps (1-q) \, \frac{a}{a+o(a)} + \frac{c(F,M,q)}{(1+a)^{1-q}-1} + \eps q\, 
   \frac{\lambda_q([M,a])}{\lambda_q([0,a])} ,
\end{split}
\end{equation}
as $a \to \infty$. Here $c(F,M,q)$ is a constant that depends on $F$, $M$ 
and $q$, but not on $a$. It follows that
\begin{equation} 
  \limsup_{a\to \infty} \left| \frac{1}{\lambda_q([0,a])} \int_0^a F\, d\lambda_q 
  \right| \le \eps (1-q) + \eps q = \eps ,
\end{equation}
for all $\eps>0$, proving (\ref{voglio-q}).

In the case $q=1$, we use that $\lambda_1([0,a]) = \log(1+a)$. From 
(\ref{interm-q}) we derive
\begin{equation} 
\begin{split}
  \left| \frac{1}{\lambda_1([0,a])} \int_0^a F\, d\lambda_1 \right| &\le \eps\, 
  \frac{a}{(1+a) \log(1+a)} \\
  &\quad + \frac{c(F,M)}{\log(1+a)} \\
  &\quad + \eps\, \frac{\log(1+a)-\log(1+M)}{\log(1+a)} ,
\end{split}
\end{equation}
with the obvious meaning of $c(F,M)$. The limit (\ref{voglio-q}) with $q=1$ follows.
\qed

To complete the proof of Lemma \ref{lem-claudio}\emph{(a)} we need the 
following
\begin{lemma} \label{lem-2-3-app-c}
  The following statements hold true:
  \begin{enumerate}
    \item[(i)] For $q\in (0,\frac 12)$, $\go(\rp,\lambda_q)\subseteq 
    \go(\rp,\leb)$.
    \item[(ii)] For any fixed $q_1\in (0,1)$, $\go(\rp,\lambda_{q_{_2}})
    \subseteq \go(\rp,\lambda_{q_{_1}})$ for all $q_2\in (q_1,\frac{1+q_1}{2})$.
  \end{enumerate}
\end{lemma}

\proof We first remark that \emph{(i)} can be thought of as a particular 
case of \emph{(ii)} with $q_1=0$. For the sake of clarity we prefer to state 
the two results separately.

Let us first prove \emph{(i)}. We show that $F\not\in \go(\rp,\leb)$ implies 
$F\not\in \go(\rp,\lambda_q)$. For $F\in L^\infty(\rp)$, the statement  
$F\not\in \go(\rp,\leb)$ is equivalent to saying that 
\begin{equation} \label{maj1}
  A := \limsup_{a\to \infty} \frac1a \int_0^a F(y)\, dy > B := \liminf_{a\to \infty}
  \frac1a \int_0^a F(y)\, dy .
\end{equation}
Thus we can find two increasing, diverging sequences $( \sigma_k )_{k\in\N}$ 
and $( \tau_k )_{k\in\N}$ in $\rp$ such that 
\begin{equation} \label{subsq}
  \lim_{k\to \infty} \frac{1}{\sigma_k} \int_0^{\sigma_k} F(y)\, dy = A ,
  \qquad \lim_{k\to \infty} \frac{1}{\tau_k} \int_0^{\tau_k} F(y)\, dy = B .
\end{equation}
Using (\ref{serve-dopo}) with $a=\sigma_k$ we obtain
\begin{equation} \label{marco5}
\begin{split}
  & \frac{1}{\lambda_q([0,\sigma_k])} \int_0^{\sigma_k} F(y)\, 
  \frac{1}{(1+y)^q}\, dy \\[4pt]
  &\qquad = \frac{(1-q) (1+\sigma_k)^{-q} \sigma_k}
  {(1+\sigma_k)^{1-q}-1}\, \frac{1}{\sigma_k} \int_0^{\sigma_k} F(y)\, dy \\
  &\qquad\quad + \frac{q}{\lambda_q([0,\sigma_k])} \int_0^{\sigma_k} 
  \frac{1}{(1+y)^{q+1}} \,\Big(\int_0^y F(s)\, ds \Big) dy.
\end{split}
\end{equation}

Now fix $\eps>0$. By the first limit in (\ref{subsq}), there exists $\bar k \in \N$
such that, for all $k \ge \bar k$,
\begin{equation}
  \frac 1 {\sigma_k} \int_0^{\sigma_k} F(y)\, dy > A-\eps .
\end{equation}
Also, by the definition of $B$ in (\ref{maj1}), there exists $\bar y \in \R^+$ 
such that, for all $y>\bar y$,
\begin{equation}
  \frac 1 y \int_0^y F(s)\, ds > B-\eps.
\end{equation}
Since $(\sigma_k)$ is diverging, we find $K \ge \bar k$ such that
$\sigma_k > \bar y$, for all $k \ge K$. Therefore, by (\ref{marco5}),
\begin{equation}  \label{ora}
\begin{split}
  & \frac{1}{\lambda_q([0,\sigma_k])} \int_0^{\sigma_k} 
  \frac{F(y)}{(1+y)^q}\, dy \\[4pt]
  &\qquad > \frac{(1-q) (1+\sigma_k)^{-q} \sigma_k}
  {(1+\sigma_k)^{1-q}-1}\, (A-\eps) + \frac{q}{\lambda_q([0,\sigma_k])}
  \int_0^{\sigma_k} \frac{(B-\eps) y}{(1+y)^{q+1}}\, dy ,
  \end{split}
\end{equation}
for all $k \ge K$. Rewriting $(B-\eps) y = (B-\eps)( (1+y) - 1)$, we see that  
\begin{equation} \label{marco10}
\begin{split}
  &\frac{q}{\lambda_q([0,\sigma_k])}
  \int_0^{\sigma_k} \frac{(B-\eps) y}{(1+y)^{q+1}}\, dy \\
  &\qquad = q (B-\eps)\, \frac{\lambda_q([0,\sigma_k]) 
  - \int_0^{\sigma_k} (1+y)^{-q-1}\, dy}
  {\lambda_q([0,\sigma_k])}
\end{split} 
\end{equation} 
which, together with (\ref{ora}), leads to
\begin{equation} 
  \liminf_{k\to \infty} \frac{1}{\lambda_q([0,\sigma_k])} \int_0^{\sigma_k} 
  \frac{F(y)}{(1+y)^q}\, dy > A(1-q) + Bq -\eps.
\end{equation}
Since the above l.h.s.\ does not depend on $\eps$, we conclude that
\begin{equation}
 \liminf_{k\to \infty} \frac{1}{\lambda_q([0,\sigma_k])} \int_0^{\sigma_k} 
  \frac{F(y)}{(1+y)^q}\, dy \ge A(1-q) + Bq.
\end{equation}

Analogously, putting $a=\tau_k$ in (\ref{serve-dopo}) and writing estimates 
from the above that are specular to the ones used above, we obtain
\begin{equation}
  \limsup_{k\to \infty} \frac{1}{\lambda_q([0,\tau_k])} \int_0^{\tau_k} 
  \frac{F(y)}{(1+y)^q}\, dy \le B (1-q) + A q .
\end{equation}
Since $A>B$ and $q<\frac 12$, we have that $A (1-q) + B q > 
B (1-q) + A q$. We have therefore produced two subsequences of 
$a \mapsto \lambda_q([0,a])^{-1} \int_0^a\, F\, d\lambda_q$ with 
different limits. This shows that $F\not\in \go(\rp,\lambda_q)$ and the 
proof of \emph{(i)} is finished.

\skippar

We now show \emph{(ii)} using a similar reasoning. Let $F \in L^\infty(\rp)$ 
with $F\not\in \go(\rp,\lambda_{q_{_1}})$. This means that
\begin{equation} \label{maj2}
  A := \limsup_{a\to \infty} \frac1 {\lambda_{q_{_1}}([0,a])} \int_0^a F\, 
  d\lambda_{q_{_1}} > B := \liminf_{a\to \infty} \frac1 
  {\lambda_{q_{_1}}([0,a])} \int_0^a F\, d\lambda_{q_{_1}},
\end{equation}
implying the existence of  two increasing sequences $( \sigma_k )_{k\in\N},
( \tau_k )_{k\in\N} \subset \rp$ such that $\sigma_k,\tau_k \to \infty$ and
\begin{equation}
  \lim_{k\to \infty} \frac{1}{\lambda_{q_{_1}}([0,\sigma_k])} 
  \int_0^{\sigma_k} F \, d\lambda_{q_{_1}} = A , \qquad 
  \lim_{k\to \infty} \frac{1}{\lambda_{q_{_1}}([0,\tau_k])} \int_0^{\tau_k} F 
  \, d\lambda_{q_{_1}} = B .
\end{equation}
We need an adapted version of (\ref{serve-dopo}). Integrating by parts 
we have
\begin{equation}
\begin{split}
  \int_0^a F\, d\lambda_{q_{_2}} &= \int_0^a F(y)\, \frac{1}{(1+y)^{q_2}}\, dy \\
  &= (1+a)^{q_1-q_2} \int_0^a \frac{F(y)}{(1+y)^{q_1}}\, dy \\
  &\quad+ (q_2-q_1) \int_0^a \frac{1}{(1+y)^{q_2-q_1+1}}\, \Big( \int_0^y 
  \frac{F(s)}{(1+s)^{q_1}}\, ds \Big) dy ,
\end{split}
\end{equation}
where we recall that $q_2>q_1$. By means of (\ref{lambdaq-0a}), we write 
\begin{equation} \label{serve-qui}
\begin{split}
  &\frac{1}{\lambda_{q_{_2}}([0,a])} \int_0^a F\, d\lambda_{q_{_2}} = \\
  &\qquad = \frac{1-q_2}{1-q_1} \, \frac{(1+a)^{1-q_2}-(1+a)^{q_1-q_2}}
  {(1+a)^{1-q_2}-1}\, \frac{1}{\lambda_{q_{_1}}([0,a])} \int_0^a F(y)\, 
  d\lambda_{q_{_1}}(y) \\
  &\qquad \quad + \frac{q_2-q_1}{\lambda_{q_{_2}}([0,a])} \int_0^a 
  \frac{1}{(1+y)^{q_2-q_1+1}}\, \Big( \int_0^y \frac{F(s)}{(1+s)^{q_1}}\, ds 
  \Big) dy .
\end{split}
\end{equation}
Now set $a=\sigma_k$ in (\ref{serve-qui}). Starting from (\ref{maj2})
one can repeat the same arguments that were used to show \eqref{ora} 
to prove that, for all $\eps>0$,
\begin{equation} \label{marco15}
\begin{split}
  & \liminf_{k\to \infty} \frac{1}{\lambda_{q_{_2}}([0,\sigma_k])} 
  \int_0^{\sigma_k} F(y)\, d\lambda_{q_{_2}}(y) \\
  &\qquad > (A-\eps)\, \frac{1-q_2}{1-q_1} + \liminf_{k\to \infty} \frac{q_2-q_1}
  {\lambda_{q_{_2}}([0,\sigma_k])} \int_0^{\sigma_k} \frac{(B-\eps)\, 
  \lambda_{q_{_1}}([0,y])} 
  {(1+y)^{q_2-q_1+1}}\, dy .
\end{split}
\end{equation}
In analogy with (\ref{marco10}), and using (\ref{lambdaq-0a}), we conclude 
that
\begin{equation} 
\begin{split}
  &\liminf_{k\to \infty} \frac{1} {\lambda_{q_{_2}}([0,\sigma_k])} 
  \int_0^{\sigma_k} \frac{(B-\eps)\, \lambda_{q_{_1}}([0,y])} 
  {(1+y)^{q_2-q_1+1}}\, dy \\
  &\qquad = (B-\eps)\, \liminf_{k\to \infty}  \frac{ \int_0^{\sigma_k} (1+y)^{-q_2}\, 
  dy - \int_0^{\sigma_k} (1+y)^{-q_2+q_1-1}\, dy} {(1-q_1) \, 
  \lambda_{q_{_2}} ([0,\sigma_k])} \\
  &\qquad = \frac{B-\eps} {1-q_1} .
\end{split} 
\end{equation} 
Once again, using that the l.h.s.\ of (\ref{marco15}) does not depend on $\eps$, 
we obtain
\begin{equation}
  \liminf_{k\to \infty} \frac{1}{\lambda_{q_{_2}}([0,\sigma_k])} \int_0^{\sigma_k}
  F(y)\, d\lambda_{q_{_2}}(y) \ge \frac{ A (1-q_2) + B (q_2-q_1) }{1-q_1} .
\end{equation}

Just as before, one can set $a=\tau_k$ in (\ref{serve-qui}) and take the 
limsup of the two sides to find
\begin{equation}
  \limsup_{k\to \infty} \frac{1}{\lambda_{q_{_2}}([0,\tau_k])} \int_0^{\tau_k} 
  F(y)\, d\lambda_{q_{_2}}(y) \le \frac{ B (1-q_2) + A (q_2-q_1) }{1-q_1} .
\end{equation}
On the other hand, for $q_1 \in (0,1)$ and $q_2\in (q_1,\frac{1+q_1}{2})$, 
the inequality 
\begin{equation}
   \frac{ A (1-q_2) + B (q_2-q_1) }{1-q_1}  >
   \frac{ B (1-q_2) + A (q_2-q_1) }{1-q_1} 
\end{equation}
is equivalent to $A>B$. So, in analogy to what was done earlier, we have 
found two subsequences of 
$a \mapsto \lambda_{q_{_2}}([0,a])^{-1} \int_0^a\, F\, d\lambda_{q_{_2}}$
with different limits. Then $F \not\in  \go(\rp,\lambda_{q_{_2}})$, ending the 
proof of \emph{(ii)}.
\qed

\medskip

\proofof{Lemma \ref{lem-claudio}\emph{(a)}} Let $q\in (0,1)$. The inclusion 
$\go(\rp,\leb)\subseteq \go(\rp,\lambda_q)$ and the equality 
$\overline{\lambda}_q(F) = \avgleb(F)$, for all $F \in \go(\rp, \leb)$, were
proved in Proposition \ref{prop-1-app-c}. 

As for the opposite inclusion, it was proved in Lemma 
\ref{lem-2-3-app-c}\emph{(i)}, for the case $q\in (0,\frac 12)$.
When $q=\frac 12$, using Lemma \ref{lem-2-3-app-c}\emph{(ii)} with 
$q_1 := \frac 14$ and $q_2 := q=\frac 12 \in (\frac 14, \frac 58)$, we get
\begin{equation}
  \go(\rp,\lambda_{1/2}) \subseteq \go(\rp,\lambda_{1/4}) \subseteq \go(\rp,\leb).
\end{equation}
Suppose instead that $q\in (\frac 12,\frac 34)$. By Lemma 
\ref{lem-2-3-app-c}\emph{(ii)} with $q_1 := \frac 12$ and $q_2 := q \in 
(q_1,\frac{1+q_1}{2})$, we have
\begin{equation}
  \go(\rp,\lambda_{q}) \subseteq \go(\rp,\lambda_{1/2}) \subseteq \go(\rp,\leb).
\end{equation}
We can now iterate the argument to obtain 
$\go(\rp, \lambda_q)\subseteq \go(\rp,\leb)$, for all $q\in [1-2^{-n},1-2^{-n-1})$,
for all $n\ge 1$. Assertion \emph{(a)} is proved. 
\qed

\medskip

\proofof{Lemma \ref{lem-claudio}\emph{(b)}} The facts that 
$\go(\rp,\leb) \subseteq \go(\rp,\lambda_1)$ and $\overline{\lambda}_1(F) = 
\avgleb(F)$, for all $F \in \go(\rp, \leb)$, are proved in Proposition 
\ref{prop-1-app-c}. It remains to show that the inclusion is actually strict, that is, 
there exists $F \in \go(\rp,\lambda_1) \setminus \go(\rp,\leb)$.

For $k \in \N$, set
\begin{equation}
  \alpha_k := k^k-1 , \qquad \beta_k := 2k^k-1 .
\end{equation}
Notice that $\alpha_1=0$, $\alpha_k < \beta_k < \alpha_{k+1}$ for all $k$,
and $\alpha_k, \beta_k \to \infty$, as $k \to \infty$. Now define a \fn\
$F : \rp \into \R$ as follows. For all $k \in \N$,
\begin{equation} \label{es-claudio}
  F(y) := \left\{ \begin{array}{ll} 
    1 , & y\in [\alpha_k,\beta_k) ; \\
    0 , & y\in [\beta_k,\alpha_{k+1}) ,
  \end{array} \right.
\end{equation}

We first show that $F\in \go(\rp,\lambda_1)$ and in particular that 
\begin{equation} \label{scopo-q1}
  \overline{\lambda}_1(F) = \lim_{a\to \infty} \frac{1}{\lambda_1([0,a])} \int_0^a 
  F(y)\, \frac{1}{1+y}\, dy = 0 .
\end{equation}
Letting $a=\alpha_n$, for $n> 1$, we find
\begin{equation}
\begin{split}
  \frac{1}{\lambda_1([0,\alpha_n])} \int_0^{\alpha_{_n}} F(y)\, \frac{1}{1+y}
  \, dy \, &= \frac{1}{\log(1+\alpha_n)} \sum_{k=1}^{n-1} 
  \int_{\alpha_k}^{\beta_k} \frac{1}{1+y}\, dy \\
  &= \frac{1}{\log(1+\alpha_n)} \sum_{k=1}^{n-1} \log \left( \frac{1+\beta_k}
  {1+\alpha_k}\right) \\
  &= \frac{1}{n \log n}\, \sum_{k=1}^{n-1} \log 2 = O \! \left( 
  \frac 1 {\log n}\right),
\end{split}
\end{equation}
as $n \to \infty$. In the same way, for $a=\beta_n$, we find
\begin{equation}
\begin{split}
  \frac{1}{\lambda_1([0,\beta_n])} \int_0^{\beta_{_n}} F(y) \frac{1}{1+y}\, dy 
  &= \frac{1}{\log(1+\beta_n)} \sum_{k=1}^{n} \int_{\alpha_k}^{\beta_k} 
  \frac{1}{1+y}\, dy \\
  &= \frac{1}{\log(1+\beta_n)} \sum_{k=1}^{n} \log \left( \frac{1+\beta_k}
  {1+\alpha_k}\right) \\
  &= \frac{1}{\log 2 + n \log n} \sum_{k=1}^{n} \log 2 = O \! \left( 
  \frac 1 {\log n}\right) .
\end{split}
\end{equation}
Moreover, 
\begin{align}
  \forall a\in (\alpha_n,\beta_n), \ \ &\frac{1}{\lambda_1([0,a])} \int_0^a 
  F(y)\, d\lambda_1(y) \le \frac{1}{\lambda_1([0,\alpha_n])} \int_0^{\beta_n}
  F(y)\, d\lambda_1(y) ; \\[4pt]
  \forall a\in (\beta_n,\alpha_{n+1}), \ \ &\frac{1}{\lambda_1([0,a])} \int_0^a
  F(y)\, d\lambda_1(y) \le \frac{1}{\lambda_1([0,\beta_n])} 
  \int_0^{\alpha_{n+1}} \!\! F(y)\, d\lambda_1(y) .
\end{align}
Since $\lambda_1([0,\alpha_n]) \sim \lambda_1([0,\alpha_{n+1}]) \sim 
\lambda_1([0,\beta_n])$, the above l.h.sides vanish, as $n \to \infty$. 
Hence \eqref{scopo-q1} is proved and $F\in \go(\rp,\lambda_1)$. 

Lastly, we prove that $F \not\in \go(\rp,\leb)$ by showing that the 
infinite-volume average
\begin{equation}
  \avgleb(F) = \lim_{a\to\infty} \, \frac 1 a \int_0^a F(y)\, dy
\end{equation}
does not exist. In fact, let $a=\alpha_n$, with $n>1$. Then
\begin{equation}
\begin{split}
  \frac{1}{\alpha_n} \int_0^{\alpha_n} F(y)\, dy &= \frac{1}{\alpha_n} 
  \sum_{k=1}^{n-1} \int_{\alpha_k}^{\beta_k}\ dy = \frac{1}{\alpha_n} 
  \sum_{k=1}^{n-1} (\beta_k-\alpha_k) \\
  &= \frac{1}{n^n-1} \sum_{k=1}^{n-1} k^k \le \frac{1}{n^n-1} 
  \sum_{k=1}^{n-1} n^k = O \! \left( \frac 1 {n}\right),
\end{split}
\end{equation}
as $n \to \infty$. On the other hand, for $a=\beta_n$ and 
$n>1$, we see that
\begin{equation}
\begin{split}
  \frac{1}{\beta_n} \int_0^{\beta_n} F(y)\, dy &= \frac{1}{\beta_n} 
  \sum_{k=1}^{n} \int_{\alpha_k}^{\beta_k} dy = \frac{1}{\beta_n} 
  \sum_{k=1}^{n} (\beta_k-\alpha_k) \\
  &= \frac{1}{2n^n-1} \sum_{k=1}^{n}\, k^k = \frac{n^n}{2n^n-1} + 
  \frac{1}{2n^n-1} \sum_{k=1}^{n-1} k^k \ge \frac 12 .
\end{split}
\end{equation}
Therefore $F \not\in \go(\rp,\leb)$ and the lemma is proved. 
\qed

\end{document}